\documentclass[14t]{article}
\usepackage{amsmath,amssymb,amsthm}
\usepackage{times}
\newtheorem{definition}{Definition}
\newtheorem{remark}[definition]{Remark}
\newtheorem{lemma}[definition]{Lemma}

\newtheorem{proposition}[definition]{Proposition}

\begin{document}
\global\def\refname{{\normalsize \it References:}}
\baselineskip 12.5pt
%
%
%
\title{\LARGE \bf Deformations and elements of deformation theory}

\date{}

\author{\hspace*{-10pt}
\begin{minipage}[t]{2.7in} \normalsize \baselineskip 12.5pt
\centerline{NIKOLAJ GLAZUNOV}
\centerline{National Aviation University}
\centerline{Department of Electronics}
\centerline{1 Cosmonaut Komarov Avenue, 03680 Kiev}
\centerline{UKRAINE}
\centerline{glanm@yahoo.com}
\end{minipage} \kern 0in
%
%
\\ \\ \hspace*{-10pt}
\begin{minipage}[b]{6.9in} \normalsize
\baselineskip 12.5pt {\it Abstract:}
This article consisted of an elementary introduction to deformation theory of varieties, schemes and manifolds, with some applications to local and global shtukas and fever to Newton polygons of $p$-divisible groups .
Soft problems and results mainly are considered. In the framework we give review of some novel results in  the theory of local shtukas, Anderson-modules,  global shtukas,  Newton polygons of $p$-divisible groups and on deformations of $p$-divisible groups with given Newton polygons
\\ [4mm] {\it Key--Words:}
Dual number, infinitesimal deformation, 
Drinfeld module, local Anderson-module, local shtuka, global shtuka, moduli stack, formal Lie group, cotangent complex, commutative group scheme,
uniformization, rigidity,  Newton polygon, loop group of a reductive group 
\end{minipage}
\vspace{-10pt}}

\maketitle

\thispagestyle{empty} \pagestyle{empty}
%
%
\section{Introduction}
\label{S1} \vspace{-4pt}

This article consisted of an elementary introduction to deformation theory of varieties, schemes and manifolds, with some applications to local and global shtukas and fever to Newton polygons of $p$-divisible groups .
From the point of view of rigid, hard and soft problems and results in this paper we consider  mainly soft problems and results. In the framework we review  some novel results and methods in the theory of local shtukas, Anderson-modules,  global shtukas and  Newton polygons of $p$-divisible groups.
 These include (but not exhaust) 
methods and results by V.~Drinfeld\cite{Dr},  U.~Hartl, E.~Viehmann \cite{HV}, R.~Singh\cite{Singh}, A.~Rad \cite{Rad}, S.~Harashita\cite{Hara} and others.
M. Gromov\cite{Gr} in his talk at the International Congress of Mathematicians in Berkeley have presented problems and results of soft and hard symplectic geometry. In this connection  
we will present some soft problems and results in dynamics and in arithmetic geometry. "Soft"  problems and results  in our talk are limited to the framework of deformations, infinitesimal deformations,  and elements of local Anderson-modules, local shtukas and global shtukas.
Review of some novel results and methods on rigidity in arithmetic geometry and in dynamics is given in author's papers  \cite{Gl1,Gl2}.

\section{Infinitesimal neighborhoods  and  infinitesimal deformations}
\label{S2} \vspace{-4pt}
Dual numbers were introduced and were used in W. Clifford~\cite{Cl} ,  E. Study~\cite{St}, R. von-Mises~\cite{vMi}.
Expand following to  \cite{Hartshorne,Shaph} the definition of dual numbers so that it was true over any field .
Let $k$ be a field, and $k[\epsilon]$ - the ring of polynomials over $k$ in the variable $\epsilon$.
Factoring $k[\epsilon]$ by the ideal $(\epsilon^ 2)$ to give the desired ring of dual numbers $D$ over $k: D = k[\epsilon] / (\epsilon^ 2)$.
Like the classical ring of dual numbers, ring $D$ is a nilpotent ring.
Further, unless specifically stated, we are under dual numbers mean elements of the above-defined rings $D$.
Let $Spec \, D$ be the corresponding affine scheme. The scheme has one geometric point which corresponds to the maximal ideal $(\epsilon)$. Its structure sheaf has nilpotent elements that distinguishes the scheme from classical algebraic varieties.

\subsection{Infinitesimal neighborhoods}
Let $X$ be a scheme over an algebraically closed field $k.$ 
Here and in the following subsection we mean  under a point of a scheme its geometric (closed) point, unless otherwise stated.
We fix, following to \cite{Shaph}], notations: $o$ - closed point of the scheme $Spec \, k$, $\overline o$ - closed point scheme $Spec \, D$, $i: Spec \, k \to Spec \, D$ - canonical embedding, under which ${\overline o} = i(o)$.
\begin{lemma}
Any morphism $\phi: Spec \, D \to X$ defines a morphism $\phi \circ i(o): Spec \, k \to X$, where $x = \phi \circ i(o)$ is a closed point of $X.$
\end{lemma}
{\em Proof}. The homomorphism $D \to k$ with kernel $(\epsilon)$ defines the canonical embedding $i$, and the morphism $\phi \circ i(o)$ defines a closed point $x \in X.$ 

Let $U$ be an affine neighborhood of the point $x, \, {\mathfrak m}_x$ be the maximal ideal of the point $x$ in $k[U]$. 

\begin{lemma}
Let ${\mathfrak M}_x $ be the set of morphisms of schemes $ Spec \, D \to X$  such that $\phi ({\overline o}) = x$. Then ${\mathfrak M}_x(Spec \, D, X) = {\mathfrak M}_x(Spec \, D, U)$.
\end{lemma}
The proof follows from local properties of schemes.

\begin{definition}
The scheme $T_x = Spec \, k[U]/ {\mathfrak m}_x^2$ is called the infinitesimal neighborhood of  the first order of the point $x$.
\end{definition}

\begin{remark}
The homomorphism $ Spec \, k[U] \to Spec \, k[U]/ {\mathfrak m}_x^2$  defines the closure imbedding $T_x \to U$ and  $T_x$ is the closed subscheme in $U$.
\end{remark}

Using these statements we have
\begin{proposition}
Morphisms $ Spec \, D \to X$  that transform $Spec \, k$ in $x \in X$ are in one-to-one correspondence with morphisms $ Spec \, D \to T_x$.
\end{proposition}

\subsection{Infinitesimal deformations}
Let $X \to Y$ be a morphism of schemes. A scheme $X$ is flat over $Y$ if the sheaf ${\cal O}_X$ is flat over ${\cal O}_Y$ \cite{Hartshorne,Shaph}. 

\begin{definition}
Let $X, T$ be schemes and $X \to T$ the flat morphism with fixed point $t \in T$ such that $X_t \simeq X_o$. In the conditions the scheme $X$ is called the (global) deformation of the scheme   $ X_o$.
\end{definition}

\begin{definition}
Let $X_0$ be the scheme of finite type over field $k$ and $D$ the ring of dual numbers over $k.$ In the conditions the scheme $X^{'}$ which is flat over  $D$ and such that $X^{'} {\otimes_D} k  \simeq X_o $  is called the infinitesimal deformation of the scheme   $ X_o$.
\end{definition}

\begin{proposition}
Given a global deformation of the scheme $ X_o$ ,  then there exists an infinitesimal deformation of the scheme   $ X_o$.
\end{proposition}
{\em Proof}. By Lemma 1 there exists the morphism $ Spec \, D \to T$  which is defined by some element of the tangent space to $T$ in point $t$.
Hence there is a scheme $X^{'}$  flat over $D$ with the closed fiber  $ X_o$ so that $X^{'} \otimes_D k  \simeq X_o $. It gives the required infinitesimal deformation.

\section{On Local shtukas and divisible local Anderson-modules}
\label{S4} \vspace{-4pt}

In their paper~\cite{HV} U. Hartl, E. Viehmann have  investigated deformations and moduli spaces of bounded local $G-$shtukas. Latest (bounded local $G-$shtukas) are function field analogs for $p-$divisible groups with extra structure.

The author~\cite{Singh} investigates relation between finite shtukas and strict finite flat commutative group schemes and relation between divisible
 local Anderson modules and formal Lie groups. Let $Nilp_{{\mathbb F}_{q}[[\xi]]}$ be the category of ${\mathbb F}_{q}[[\xi]]$-schemes on which $\xi$ is locally nilpotent. Let $S \in Nilp_{{\mathbb F}_{q}[[\xi]]}$. The main result of this dissertation by R.~Singh is the following (section 2.5) interesting result: it is possible to associate a formal Lie group to any $z$-divisible local Anderson module over $S$ in the case when $\xi$ is locally nilpotent on $S$.


A general framework for the dissertation is the decent theory by A. Grothendieck~\cite{Gro}.  and his colleagues,  its extensions and  specializations to finite characteristics.

In Chapter 1 the author of the dissertation~\cite{Singh} defines cotangent complexes as in papers by S. Lichtenbaum and M. Schlessinger~\cite{LiSch} , by W. Messing~\cite{Me} , by V. Abrashkin~\cite{Ab}  and prove that they are homotopically equivalent.

More generally to any morphism $f: A \to B$ of commutative ring objects in a topos is associated a cotangent complex $L_{B/A}$ and to any morphism of commutative ring objects in a topos of finite and locally free $Spec (A)$-group scheme $G$ is associated a cotangent complex $L_{G/Spec (A)}$ as has presented in books by L. Illusie~\cite{Ill} .

In section 1.5 the author of the dissertation~\cite{Singh} investigates the deformations of affine group schemes follow to the mentioned paper of Abrashkin and defines strict finite $O-$module schemes.
Next section is devoted to relation between finite shtukas by V. Drinfeld~\cite{Dr}  and strict finite flat commutative group schemes.
The comparison between cotangent complex and Frobenius map of finite ${\mathbb F}_p$-shtukas is given in section 1.7.

$z-$divisible local Anderson modules by U. Hartl~\cite{Hart}  and local schtukas are investigated in Chapter 2.

Sections 2.1, 2.2 and 2.3 on formal Lie groups, local shtukas and divisible  local Anderson-modules define and illustrate notions for later use. Many of these, if not new, are set in a new form, 

In Section 2.4 the equivalence between the category of effective local shtukas over $S$ and the category of $z$-divisible local Anderson modules over $S$ is treated.

In the last section the theorem about canonical ${\mathbb F}_{q}[[z]]$-isomorphism of $z$-adic Tate-module of $z$-divisible local Anderson module $G$ of rank $r$ over $S$ and Tate module of local shtuka over $S$ associated to $G$ is given.

\section{On uniformizing the moduli stacks of global $\mathfrak G$-shtukas}
\label{S5} \vspace{-4pt}

The dissertation by Arasteh Rad ~\cite{Rad} is a Ph.D. Thesis, written under U. Hartl (M\"unster). The dissertation is devoted to the development of the theory of local $\mathbb{P}$-shtukas with the aim of their relation to the moduli stack of global ${\mathfrak G}$-shtukas. Here  ${\mathbb P}$ is a paraholic Bruhat-Tits group scheme by Pappas, Rapoport~\cite{PaRa}  and ${\mathfrak G}$ is a parahoric Bruhat-Tits group scheme over a smooth projective curve over finite field ${\mathbb F}_q$ with $q$ elements of characteristic $p.$ 

Let $C$ be a smooth projective geometrically irreducible curve over ${\mathbb F}_q$. A global ${\mathfrak G}$-shtuka ${\overline{\cal G}}$ over an ${\mathbb F}_q$-scheme $S$ is a tuple $({\cal G}, s_1, \ldots , s_n, \tau)$ consisting of a ${\mathfrak G}$-torsor ${\cal G} $ over 
$C_S := C {\times}_{{\mathbb F}_q}  S  $, an $n$-tuple of (characteristic) sections $( s_1, \ldots , s_n) \in C^n(S)$ and a Frobenius connection $\tau $ defined outside the graphs of the sections $s_i$ by Hartl, Rad~\cite{HR}.

Local $G$-shtukas by Hartl, Viehmann~\cite{HV}  and by Viehmann~\cite{Ve}  are generalizations to arbitrary reductive groups of the local analogue of Drinfeld shtukas.

Drinfeld Shtukas (the space $FSh_{D,r}$ of $F$-sheaves) was considered by Drinfeld~\cite{Dr} and by Lafforgue~\cite{Laff} .

For more results concerning local shtukas and Anderson-modules see dissertation by Singh~\cite{Singh}  written also under U. Hartl.

The main results of the dissertation \cite{Rad} are the following. The analogue of the Serre-Tate theorem over function fields that relating the deformation theory of global ${\mathfrak G}$-shtukas to the deformation theory of the associated local ${\mathbb P}_{\nu}$-shtukas via the global-local functor  (Theorem 4.1). Representablity of the Rappoport-Zink functor (Theorem 6.3.1.). The uniformization theorem from Section 7. Finally, the discussion about uniformization and local model of the moduli of global ${\mathfrak G}$-shtukas are given.

\section{On the supremum of Newton polygons of $p$-divisible groups with a given $p$-kernel type}
\label{S6} \vspace{-4pt}

The author of the paper~\cite{Hara} proves the existence of the supremum of Newton polygons of $p$-divisible groups with a given $p$-kernel type and provides an algorithm determining it.
The main results of the paper~\cite{Hara} are the following Theorem 1.1. $\xi(w)$ is the biggest one of the Newton polygons $\xi$ with $\mu(\xi) \subset w$., and Corollary 2.2. There exists the supremum of Newton polygons of $p$-divisible groups with the given $p^m$-kernel type.

Let $k$ be an algebraically closed field of characteristic $p > 0$, $c$ and $d$ be non-negative integers with $r := c + d > 0$. Let $W$ be the Weyl group of the general linear group $GL_r$, $s_i \in W$ the simple reflection, $ S = \{s_1, \ldots s_{r - 1}  \}$, $J := S \setminus \{ s_{\alpha}\}$ and let $w$ be any element of the set $(J, \emptyset)$-reduced elements of $W$ by 
N. Bourbaki~\cite{Bu} . 
The theorem  is an unpolarized  analogue of Corollary II by Harashita~\cite{Hara2}.

In the polarized case, the existence of the supremum $\xi(w)$ follows from the results by Ekedahl and van der Geer \cite{EG}.
In the case there is a good moduli space $A_g$ of principle polarized abelian varieties. In the unpolarized case there is no a good moduli space like $A_g$.

The difference of the author method in comparison with the Ekedahl - van der Geer approach is the using of ${\bf T}_{m}$-action by Vasiu~\cite{Va} which gives that the set of ${\bf T}_{m}$-orbits is naturally bijective to the set of isomorphism classes of truncated Barsotti-Tate groups of level $m$ over $k$ with codimension $c$ and dimension $d$.

\section{On the Newton strata in the loop group of a reductive group}
\label{S7} \vspace{-4pt}

The author of the paper~\cite{Ve}  generalizes purity of the Newton stratification to purity for a single break point of the Newton 
point in the context of local $G$-shtukas respectively of elements of the loop group of a reductive group. 
As an application she proves that elements of the loop group bounded by a given dominant coweight satisfy a 
generalization of Grothendieck`s conjecture on deformations of $p$-divisible groups with given Newton polygons.

Let $G$ be a split connected reductive group over ${\bf F}_{p}$, let $T$ be a split maximal torus of $G$ and 
let $LG$ be the loop group of $G$ by Faltings~\cite{Fa} . 

Let $R$ be a ${\bf F}_q$-algebra and $K$ be the sub-group scheme of $LG$ with $K(R) = G(R[[z]])$. Let 
$\sigma$ be the Frobenius of $k$ over ${\bf F}_q$ and also of $k((z))$ over ${\bf F}_q ((z))$. For algebraically 
closed $k$, the set of $\sigma$-conjugacy classes $[b] = \{g^{-1}b \sigma(g) | g \in G(k((z))) \} $ of elements 
$b \in LG(k)$ is classified by two invariants, the Kottwitz point $\kappa_{G}(b)$ and the Newton point $\nu $.

The author of the paper~\cite{Ve} proves the following two main results. Theorem 1: Let $S$ be an integral
 locally noetherian scheme and let $b \in LG(S)$. Let $j \in J(\nu)$ be a break point of the Newton point $\nu$
 of $b$ at the generic point of $S$. Let $U_j$ be the open subscheme of $S$ defined by the condition that a point
 $x$  of $S$ lies in $ U_j $ if and only if ${\mathrm pr}_{(j)}(\nu_b(x)) = {\mathrm pr}_{(j)}(\nu) $. Then $ U_j $ is an affine
 $S$-scheme. 

Theorem 2: Let $\mu_1 \preceq \mu_2  \in X_{*}(T)$ be dominant coweights. Let 
$S_{\mu_1, \mu_2} = {\bigcup_{\mu_1 \preceq \mu^{`} \preceq \mu_2}}K z^{\mu^{`}}K$.  Let $[b]$ be a 
$\sigma$-conjugacy class with $\kappa_{G}(b) = \overline{\mu}_1 = \overline{\mu}_2  $ as elements of 
$\pi_{1}(G)$ and with $\nu_{b} \preceq \mu_2 $. Then the Newton stratum $N_{b} = [b] \cap S_{\mu_1, \mu_2}$ 
is non-empty and pure of codimension $\langle \rho, \mu_2 - \nu_{b} \rangle + \frac{1}{2} def(b)$ in 
$ S_{\mu_1, \mu_2}$. The closure of $ N_{b}$ is the union of all $ N_{b^`}$ for $[b^`]$ with 
$\kappa_{G}(b^{`}) = \overline{\mu}_1$ and $\nu_{b^{`}} < \nu_{b}$.

Here $\rho$ is the half-sum of the positive roots of $G$ and the defect $def(b)$ is defined as 
$\mathrm{rk}G - \mathrm{rk}_{{\bf F}_q} J_b$ where $ J_b $ is the reductive group over ${\bf F}_q$ with 
$ J_b (k((z))) = \{ g \in LG(\overline{k}) | gb = b \sigma (g) \}$ for every field $k$ containing ${\bf F}_q$ 
and with algebraically closed $\overline{k}$. 

The proof of Theorem 1 is based on a generalization of some results by Vasiu~\cite{Va2}.

An interesting feature of her method in the prove of Theorem 2 is the using of various results on the Newton 
stratification on loop groups as Theorem 1 and the dimension formula for affine Deligne-Lusztig varieties by 
G{\H o}rtz, Haines, Kottwitz, Reuman \cite{GHKR} together 
with results on lengths of chains of Newton points by Chai~\cite{Cha}. 

\section{Conclusion}
\label{S8} \vspace{-4pt}

 Deformations and elements of deformation theory of manifolds, varieties and schemes have presented.
 In the framework we have reviewed  some novel results and methods in the theory of local shtukas, Anderson-modules,  global shtukas,  Newton polygons of $p$-divisible groups and on deformations of $p$-divisible groups with given Newton polygons.
  In this connection  we have  presented some soft problems and results in dynamics and in arithmetic geometry. "Soft"  problems and results  in our considerations are limited to the framework of deformations, infinitesimal deformations,  elements of local Anderson-modules, local shtukas, global shtukas  and deformations of $p$-divisible groups with given Newton polygons.

\end{document}